\documentclass[12pt,reqno]{amsart} 
\usepackage{amssymb,amscd,url}
\usepackage{upref}     %% upref forces cross-references with \ref to roman (upright)

\begin{document}

\allowdisplaybreaks

%%%%%%%%%%%%%%%%%%%%%%%%%%%%%%%%%%%%%%%%%%%%%%%%%%%%%%%%%%%%%%%%%%%%%%
%% Title and Author Information

\title[Dynamical Degrees, Arithmetic Degrees, and \dots]
{Dynamical Degrees, Arithmetic Degrees,
  and Canonical Heights: History, Conjectures, and Future Directions}

\date{\today}
\author[Joseph H. Silverman]{Joseph H. Silverman}
\email{joseph\_silverman@brown.edu}
\address{Mathematics Department, Box 1917
         Brown University, Providence, RI 02912 USA}


%%%%%%%%%%%%%%%%%%%%%%%%%%%%%%%%%%%%%%%%%%%%%%%%%%%%%%%%%%%%%%%%%%%%%%

% \allowdisplaybreaks

\hyphenation{ca-non-i-cal semi-abel-ian}

%%%%%%%%%%%%%%%%%%%%%%%%%%%%%%%%%%%%%%%%%%%%%%%%%%%%%%%%%%%%%%%%%%%%%%
% Theorem environments

\newtheorem{theorem}{Theorem}
\newtheorem{lemma}[theorem]{Lemma}
\newtheorem{sublemma}[theorem]{Sublemma}
\newtheorem{conjecture}[theorem]{Conjecture}
\newtheorem{proposition}[theorem]{Proposition}
\newtheorem{corollary}[theorem]{Corollary}
\newtheorem*{claim}{Claim}

\theoremstyle{definition}
% The * surpresses numbering
\newtheorem{definition}[theorem]{Definition}
\newtheorem{intuition}[theorem]{Intuition}
\newtheorem*{moral}{Moral}
\newtheorem{example}[theorem]{Example}
\newtheorem{remark}[theorem]{Remark}
\newtheorem{question}[theorem]{Question}
\newtheorem{update}{2023/24 Update}

\theoremstyle{remark}
\newtheorem*{acknowledgement}{Acknowledgements}

%%%%%%%%%%%%%%%%%%%%%%%%%%%%%%%%%%%%%%%%%%%%%%%%%%%%%%%%%%%%%%%%%%%%%%

%%%%%%%% Set Up Environment for Notation %%%%%%%%%%%%%%
% This is currently set to allow quite wide items to be defined
\newenvironment{notation}[0]{%
  \begin{list}%
    {}%
    {\setlength{\itemindent}{0pt}
     \setlength{\labelwidth}{4\parindent}
     \setlength{\labelsep}{\parindent}
     \setlength{\leftmargin}{5\parindent}
     \setlength{\itemsep}{0pt}
     }%
   }%
  {\end{list}}

%%%%%%%% Set Up Environment for Parts in Theorems %%%%%%%%%%%%%%
\newenvironment{parts}[0]{%
  \begin{list}{}%
    {\setlength{\itemindent}{0pt}
     \setlength{\labelwidth}{1.5\parindent}
     \setlength{\labelsep}{.5\parindent}
     \setlength{\leftmargin}{2\parindent}
     \setlength{\itemsep}{0pt}
     }%
   }%
  {\end{list}}
% Use \Part{(a)}, instead of \item[(a)], to ensure upright font
\newcommand{\Part}[1]{\item[\upshape#1]}

%%%%%%%% Set Up Macro for Cases %%%%%%%%%%%%%%
\def\Case#1#2{%
\paragraph{\textbf{\boldmath Case #1: #2.}}\hfil\break\ignorespaces}

%%%%%%%%%%%%%%%%%%
% Greek Alphabet %
%%%%%%%%%%%%%%%%%%
\renewcommand{\a}{\alpha}
\renewcommand{\b}{\beta}
\newcommand{\g}{\gamma}
\renewcommand{\d}{\delta}
\newcommand{\e}{\epsilon}
\newcommand{\f}{\varphi}
\newcommand{\bfphi}{{\boldsymbol{\f}}}
\renewcommand{\l}{\lambda}
\renewcommand{\k}{\kappa}
\newcommand{\lhat}{\hat\lambda}
\newcommand{\m}{\mu}
\newcommand{\bfmu}{{\boldsymbol{\mu}}}
\renewcommand{\o}{\omega}
\renewcommand{\r}{\rho}
\newcommand{\rbar}{{\bar\rho}}
\newcommand{\s}{\sigma}
\newcommand{\sbar}{{\bar\sigma}}
\renewcommand{\t}{\tau}
\newcommand{\z}{\zeta}

\newcommand{\D}{\Delta}
\newcommand{\G}{\Gamma}
\newcommand{\F}{\Phi}
\renewcommand{\L}{\Lambda}

%%%%%%%%%%%%%%%%%%%%
% Fraktur Alphabet %
%%%%%%%%%%%%%%%%%%%%
\newcommand{\ga}{{\mathfrak{a}}}
\newcommand{\gb}{{\mathfrak{b}}}
\newcommand{\gn}{{\mathfrak{n}}}
\newcommand{\gp}{{\mathfrak{p}}}
\newcommand{\gP}{{\mathfrak{P}}}
\newcommand{\gq}{{\mathfrak{q}}}

%%%%%%%%%%%%%%%%%%%
% Barred Alphabet %
%%%%%%%%%%%%%%%%%%%
\newcommand{\Abar}{{\bar A}}
\newcommand{\Ebar}{{\bar E}}
\newcommand{\kbar}{{\bar k}}
\newcommand{\Kbar}{{\bar K}}
\newcommand{\Pbar}{{\bar P}}
\newcommand{\Sbar}{{\bar S}}
\newcommand{\Tbar}{{\bar T}}

%%%%%%%%%%%%%%%%%%%%%%%%%
% Calligraphic Alphabet %
%%%%%%%%%%%%%%%%%%%%%%%%%
\newcommand{\Acal}{{\mathcal A}}
\newcommand{\Bcal}{{\mathcal B}}
\newcommand{\Ccal}{{\mathcal C}}
\newcommand{\Dcal}{{\mathcal D}}
\newcommand{\Ecal}{{\mathcal E}}
\newcommand{\Fcal}{{\mathcal F}}
\newcommand{\Gcal}{{\mathcal G}}
\newcommand{\Hcal}{{\mathcal H}}
\newcommand{\Ical}{{\mathcal I}}
\newcommand{\Jcal}{{\mathcal J}}
\newcommand{\Kcal}{{\mathcal K}}
\newcommand{\Lcal}{{\mathcal L}}
\newcommand{\Mcal}{{\mathcal M}}
\newcommand{\Ncal}{{\mathcal N}}
\newcommand{\Ocal}{{\mathcal O}}
\newcommand{\Pcal}{{\mathcal P}}
\newcommand{\Qcal}{{\mathcal Q}}
\newcommand{\Rcal}{{\mathcal R}}
\newcommand{\Scal}{{\mathcal S}}
\newcommand{\Tcal}{{\mathcal T}}
\newcommand{\Ucal}{{\mathcal U}}
\newcommand{\Vcal}{{\mathcal V}}
\newcommand{\Wcal}{{\mathcal W}}
\newcommand{\Xcal}{{\mathcal X}}
\newcommand{\Ycal}{{\mathcal Y}}
\newcommand{\Zcal}{{\mathcal Z}}

%%%%%%%%%%%%%%%%%%%%%%%%%%%%
% Blackboard Bold Alphabet %
%%%%%%%%%%%%%%%%%%%%%%%%%%%%
\renewcommand{\AA}{\mathbb{A}}
\newcommand{\BB}{\mathbb{B}}
\newcommand{\CC}{\mathbb{C}}
\newcommand{\FF}{\mathbb{F}}
\newcommand{\GG}{\mathbb{G}}
\newcommand{\NN}{\mathbb{N}}
\newcommand{\PP}{\mathbb{P}}
\newcommand{\QQ}{\mathbb{Q}}
\newcommand{\RR}{\mathbb{R}}
\newcommand{\ZZ}{\mathbb{Z}}

%%%%%%%%%%%%%%%%%%%%%%%%%%
% Boldface Math Alphabet %
%%%%%%%%%%%%%%%%%%%%%%%%%%
\newcommand{\bfa}{{\boldsymbol a}}
\newcommand{\bfb}{{\boldsymbol b}}
\newcommand{\bfc}{{\boldsymbol c}}
\newcommand{\bfd}{{\boldsymbol d}}
\newcommand{\bfe}{{\boldsymbol e}}
\newcommand{\bff}{{\boldsymbol f}}
\newcommand{\bfg}{{\boldsymbol g}}
\newcommand{\bfi}{{\boldsymbol i}}
\newcommand{\bfj}{{\boldsymbol j}}
\newcommand{\bfp}{{\boldsymbol p}}
\newcommand{\bfr}{{\boldsymbol r}}
\newcommand{\bfs}{{\boldsymbol s}}
\newcommand{\bft}{{\boldsymbol t}}
\newcommand{\bfu}{{\boldsymbol u}}
\newcommand{\bfv}{{\boldsymbol v}}
\newcommand{\bfw}{{\boldsymbol w}}
\newcommand{\bfx}{{\boldsymbol x}}
\newcommand{\bfy}{{\boldsymbol y}}
\newcommand{\bfz}{{\boldsymbol z}}
\newcommand{\bfA}{{\boldsymbol A}}
\newcommand{\bfF}{{\boldsymbol F}}
\newcommand{\bfB}{{\boldsymbol B}}
\newcommand{\bfD}{{\boldsymbol D}}
\newcommand{\bfG}{{\boldsymbol G}}
\newcommand{\bfI}{{\boldsymbol I}}
\newcommand{\bfM}{{\boldsymbol M}}
\newcommand{\bfP}{{\boldsymbol P}}
\newcommand{\bfzero}{{\boldsymbol{0}}}
\newcommand{\bfone}{{\boldsymbol{1}}}

%%%%%%%%%%%%%%%%%%%%%%%%%%%%%%
% Miscellaneous New Commands %
%%%%%%%%%%%%%%%%%%%%%%%%%%%%%%
\newcommand{\Aut}{\operatorname{Aut}}
\newcommand{\codim}{\operatorname{codim}}
\newcommand{\Crit}{\operatorname{Crit}}
\newcommand{\Disc}{\operatorname{Disc}}
\newcommand{\Div}{\operatorname{Div}}
\newcommand{\Dom}{\operatorname{Dom}}
\newcommand{\dyn}{{\textup{dyn}}}
\newcommand{\End}{\operatorname{End}}
\newcommand{\Fbar}{{\bar{F}}}
\newcommand{\Fix}{\operatorname{Fix}}
\newcommand{\Gal}{\operatorname{Gal}}
\newcommand{\GL}{\operatorname{GL}}
\newcommand{\Hom}{\operatorname{Hom}}
\newcommand{\Image}{\operatorname{Image}}
\newcommand{\Isom}{\operatorname{Isom}}
\newcommand{\hhat}{{\hat h}}
\newcommand{\Ker}{{\operatorname{ker}}}
\newcommand{\limstar}{\lim\nolimits^*}
\newcommand{\limstarn}{\lim_{\hidewidth n\to\infty\hidewidth}{\!}^*{\,}}
\newcommand{\Mat}{\operatorname{Mat}}
\newcommand{\maxplus}{\operatornamewithlimits{\textup{max}^{\scriptscriptstyle+}}}
\newcommand{\MOD}[1]{~(\textup{mod}~#1)}
\newcommand{\Mor}{\operatorname{Mor}}
\newcommand{\Moduli}{\mathcal{M}}
\newcommand{\Norm}{{\operatorname{\mathsf{N}}}}
\newcommand{\notdivide}{\nmid}
\newcommand{\normalsubgroup}{\triangleleft}
\newcommand{\NS}{\operatorname{NS}}
\newcommand{\onto}{\twoheadrightarrow}
\newcommand{\ord}{\operatorname{ord}}
\newcommand{\Orbit}{\mathcal{O}}
\newcommand{\Per}{\operatorname{Per}}
\newcommand{\PrePer}{\operatorname{PrePer}}
\newcommand{\PGL}{\operatorname{PGL}}
\newcommand{\Pic}{\operatorname{Pic}}
\newcommand{\Prob}{\operatorname{Prob}}
\newcommand{\Proj}{\operatorname{Proj}}
\newcommand{\Qbar}{{\bar{\QQ}}}
\newcommand{\rank}{\operatorname{rank}}
\newcommand{\Rat}{\operatorname{Rat}}
\newcommand{\Resultant}{\operatorname{Res}}
\renewcommand{\setminus}{\smallsetminus}
\newcommand{\sgn}{\operatorname{sgn}} 
\newcommand{\SL}{\operatorname{SL}}
\newcommand{\Span}{\operatorname{Span}}
\newcommand{\Spec}{\operatorname{Spec}}
\renewcommand{\ss}{\textup{ss}}
\newcommand{\stab}{\textup{stab}}
\newcommand{\Support}{\operatorname{Supp}}
\newcommand{\tors}{{\textup{tors}}}
\newcommand{\tr}{{\textup{tr}}} 
\newcommand{\Trace}{\operatorname{Trace}}
\newcommand{\trianglebin}{\mathbin{\triangle}} % symmetric set difference
\newcommand{\UHP}{{\mathfrak{h}}}    % Upper half plane
\newcommand{\<}{\langle}
\renewcommand{\>}{\rangle}

\newcommand{\pmodintext}[1]{~\textup{(mod}~#1\textup{)}}
\newcommand{\ds}{\displaystyle}
\newcommand{\longhookrightarrow}{\lhook\joinrel\longrightarrow}
\newcommand{\longonto}{\relbar\joinrel\twoheadrightarrow}
\newcommand{\SmallMatrix}[1]{%
  \left(\begin{smallmatrix} #1 \end{smallmatrix}\right)}
%% This creates an \xrightarrow that's dashed. This is used to create a long dashed arrow.
\makeatletter
\newcommand{\xdashrightarrow}[2][]{\ext@arrow 0359\rightarrowfill@@{#1}{#2}}
\def\rightarrowfill@@{\arrowfill@@\relax\relbar\rightarrow}
\def\arrowfill@@#1#2#3#4{%
  $\m@th\thickmuskip0mu\medmuskip\thickmuskip\thinmuskip\thickmuskip
   \relax#4#1
   \xleaders\hbox{$#4#2$}\hfill
   #3$%
}
\newcommand{\longdashrightarrow}{\xdashrightarrow{\hspace{2.5em}}}
\makeatother
\def\crit{{\textup{crit}}}
\def\hhatinf{\underline{\hhat}}

%%%%%%%%%%%%%%%%%%%%%%%%%%%%%%%%%%%%%%%%%%%%%%%%%%%%%%%%%%%%%%%%%%%%%%

\begin{abstract}
In this note we give an overview of various quantities that are used
to measure the complexity of an algebraic dynamical system
$f:X\to{X}$, including the dynamical degree $\delta(f)$, which gives a
coarse measure of the geometric complexity of the iterates of $f$, the
arithmetic degree $\alpha(f,P)$, which gives a coarse measure of the
arithmetic complexity of the orbit of
$P\in{X}(\overline{\mathbb{Q}})$, and various versions of the
canonical height~$\hat{h}_f(P)$ that provide more refined measures of
arithmetic complexity. Emphasis is placed on open problems and
directions for further exploration.
\end{abstract}

\maketitle

%% \tableofcontents

This article is a slightly expanded version of a talk presented at
the Simons Symposium on Algebraic, Complex and Arithmetic Dynamics
held May 19--23, 2019 at Schloss Elmau, Kr\"un, Germany. A small
number of updates were added in~2023 and~2024 and are noted as such.

%%%%%%%%%%%%%%%%%%%%%%%%%%%%%%%%%%%%%%%%%%%%%%%%%%%%%%%%%%%%%%%%%%%%%%
\section{Measuring Complexity of Iteration}
%%%%%%%%%%%%%%%%%%%%%%%%%%%%%%%%%%%%%%%%%%%%%%%%%%%%%%%%%%%%%%%%%%%%%%

Let $X$ be an object in some category, and suppose that we are given
a function
\[
  h : \End(X) \longrightarrow \RR_{\ge0}
\]
that measures the complexity of endomorphisms of~$X$.
\begin{center}
  \framebox{\parbox{.75\hsize}{\noindent
      \textbf{The Endomorphism Complexity Problem}:\\
      Describe the growth rate of $h(f^n)$ as $n\to\infty$.
    }
  }
\end{center}

Next suppose that we are working in a category where the objects are
sets, and that for every object $X$ we are given a function
\[
  h_X : X \longrightarrow \RR_{\ge0}
\]
that measures the complexity of the elements of~$X$.

\begin{center}
  \framebox{\parbox{.75\hsize}{\noindent
      \textbf{The Orbit Complexity Problem}:\\
      For $x\in X$, describe the growth rate of $h_X\bigl(f^n(x)\bigr)$ as
      $n\to\infty$. Classify the subsets of $X$ exhibiting various growth
      rates.
    }
  }
\end{center}

This is  all very abstract, so we move on to explicit examples.

%%%%%%%%%%%%%%%%%%%%%%%%%%%%%%%%%%%%%%%%%%%%%%%%%%%%%%%%%%%%%%%%%%%%%%
\section{The Endomorphism Complexity Problem for $\PP^N$}
%%%%%%%%%%%%%%%%%%%%%%%%%%%%%%%%%%%%%%%%%%%%%%%%%%%%%%%%%%%%%%%%%%%%%%

We start with projective space. We consider
\[
  f:\PP^N\longdashrightarrow\PP^N,
  \quad\text{a dominant rational map.}
\]
One measure of complexity is the degree of a map,
\[
  \deg : \End(\PP^N) \longrightarrow \NN_{\ge1}.
\]
If $f$ is a morphism, then it is easy to show that
\[
\deg(f^n)=\deg(f)^n,
\]
but in general the growth of the degree sequence can
be quite complicated.
\begin{definition}
The \emph{dynamical degree of $f$} is the quantity
\[
  \d(f) := \lim_{n\to\infty} \Bigl(\deg(f^n)\Bigr)^{1/n}.
\]
\end{definition}

\begin{intuition}
$\deg(f^n)$ is roughly $\d(f)^n$.
\end{intuition}

\begin{conjecture}
  \label{conj:BV}
  \textup{(Bellon-Viallet \cite{MR1704282})}
  \[
  \d(f)\in\overline{\ZZ}.
  \]
\end{conjecture}

Conjecture~\ref{conj:BV} is known for various types of varieties and
maps. For example:

\begin{theorem}
  \textup{(Diller--Favre \cite{MR1867314})}
  Conjecture~$\ref{conj:BV}$ is true for birational maps $\PP^2\dashrightarrow\PP^2$.
\end{theorem}

\begin{update}
Bell, Diller, and Jonsson~\cite{arxiv1907.00675} have shown that there
are dominant rational maps $f:\PP^2\dashrightarrow\PP^2$ whose
dynamical degree~$\d(f)$ is transcendental over~$\QQ$, and Bell,
Diller, Jonsson, and Krieger~\cite{belldillerjonssonkrieger2021} have
shown that for all~$N\ge3$, there are birational maps
$f:\PP^N\dashrightarrow\PP^N$ whose dynamical degree~$\d(f)$ is
transcendental over~$\QQ$.
\end{update}

%%%%%%%%%%%%%%%%%%%%%%%%%%%%%%%%%%%%%%%%%%%%%%%%%%%%%%%%%%%%%%%%%%%%%%
\section{The Endomorphism Complexity Problem for Varieties}
%%%%%%%%%%%%%%%%%%%%%%%%%%%%%%%%%%%%%%%%%%%%%%%%%%%%%%%%%%%%%%%%%%%%%%

Let $X$ be a smooth irreducible projective variety of dimension~$N$, and let~$H$
be an ample divisor on~$X$, and let~$\End(X)$ denote the set
of dominant rational maps~$f:X\dashrightarrow{X}$. We can measure the complexity of
$f\in\End(X)$ by its $H$-degree,
\[
\deg_H : \End(X)\longrightarrow\NN_{\ge1},
\quad
\deg_H(f) := f^*H \cdot H^{N-1}.
\]

\begin{definition}
\label{definition:dynamicaldegree}
The \emph{dynamical degree of $f$} is the quantity\footnote{More
generally, one can define intermediate dynamical degrees~$\d_k(f)$
for each~$0\le{k}\le{N}$. See Section~\ref{section:higherorderdyndeg}.}
\[
\d(f) := \lim_{n\to\infty} \Bigl(\deg_H(f^n)\Bigr)^{1/n}.
\]
\end{definition}

A map~$f\in\End(X)$ induces a map on the N{\'e}ron--Severi
group~$\NS(X)$, and we adopt the usual
notation~$\NS(X)_\QQ=\NS(X)\otimes\QQ$.

\begin{itemize}
  \setlength{\itemsep}{0pt}
\item
  In general, one finds that $(f^n)^*\ne(f^*)^n$ as maps on $\NS(X)_\QQ$.
\item
  The limit $\d(f)$ exists and is independent of~$H$.
\item
  It is enough to take $X$ normal and $H$ to be a Cartier divisor that
  is nef and big; see for example~\cite[Theorem~1]{arXiv:1701.07760}.
\end{itemize}

%%%%%%%%%%%%%%%%%%%%%%%%%%%%%%%%%%%%%%%%%%%%%%%%%%%%%%%%%%%%%%%%%%%%%%
\section{Variation of the Dynamical Degree in Families}
%%%%%%%%%%%%%%%%%%%%%%%%%%%%%%%%%%%%%%%%%%%%%%%%%%%%%%%%%%%%%%%%%%%%%%

\begin{conjecture}
\label{conj:limp}
Let $f:\PP^N_\QQ\dashrightarrow\PP^N_\QQ$. For all sufficiently large primes~$p$, we can
reduce modulo~$p$ to obtain a map
\[
\tilde{f}_p:\PP^N_{\FF_p}\longdashrightarrow\PP^N_{\FF_p}.
\]
\textup(We note that $\d(\tilde f_p) \le \d(f)$ for all~$p$.\textup)
Then
\[
  \lim_{p\to\infty} \d(\tilde f_p) = \d(f).
\]
\end{conjecture}

\begin{update}
Xie~\cite[Corollary~1.10]{arxiv2402.12678} has proven Conjecture~\ref{conj:limp}.
\end{update}

\begin{remark}
Xie~\cite{arxiv1106.1825} gives an example such that there  
is a strict inequality $\d(\tilde f_p) < \d(f)$ for all $p$.
\end{remark}

More generally, Conjecture~\ref{conj:limp} should be true for the
reduction modulo~$\gp$ of dominant rational maps of smooth irreducible projective
varieties defined over number fields.  Note that we may view this
conjecture as studying the variation of a family of maps
over~$\Spec\ZZ$ as we specialize to closed fibers.  This suggests an
analogous geometric analogue.

Let $f:X/T\dashrightarrow X/T$ be a family of dominant rational maps
parameterized by an irreducible variety~$T$. This gives a dynamical degree
$\d(f_\eta)$ of $f$ on the generic fiber, i.e., over the function
field $k(T)$ of~$T$, and also, for each~$t\in T$, a dynamical degree
$\d(f_t)$ on the map $f_t:X_t\dashrightarrow{X_t}$ of the fiber over~$t$.

\begin{conjecture}
\label{conj:tinTd}
\textup{(Call--Silverman, \cite{MR3841148})}
For all $\e>0$, the set  
\[
\bigl\{t\in T : \d(f_t)\le\d(f_\eta)-\e \bigr\}
\quad\text{is not Zariski dense in $T$.}
\]
\end{conjecture}

\begin{theorem}
\textup{(Xie \cite{arxiv1106.1825})}
Conjectures~$\ref{conj:limp}$ and~$\ref{conj:tinTd}$ are true for
$\PP^2_T:=\PP^2\times{T}$, i.e., for families of dominant rational
maps~$\PP^2\dashrightarrow\PP^2$ parameterized by the points of an irreducible
variety~$T$.
\end{theorem}

\begin{update}
Let $f:X/T\dashrightarrow X/T$ be a family of dominant rational maps
parameterized by an irreducible variety~$T$ as described
in~\cite[Definition~1.7]{arxiv2402.12678}.  (Or more generally, the
parameter space~$T$ may be an integral noetherian scheme.)  Then
Xie~\cite[Theorem~1.9]{arxiv2402.12678} proves that the map
\begin{equation}
  \label{eqn:ttodeltaftlowersc}
  T(\Qbar) \longrightarrow \RR_{\ge1},\quad
  t \longmapsto \d(f_t),
\end{equation}
is lower semi-continuous with respect to the Zariski topology on~$T$
and the usual real topology on~$\RR$.\footnote{More generally, Xie
proves that~\eqref{eqn:ttodeltaftlowersc} is true for~$t\to\d_i(f_t)$
for all of the intermediate dynamical degrees as described in the
footnote to Definition~\ref{definition:dynamicaldegree}.}  A corollary
of~\eqref{eqn:ttodeltaftlowersc} is that Conjecture~\ref{conj:tinTd}
is true.
\end{update}

%%%%%%%%%%%%%%%%%%%%%%%%%%%%%%%%%%%%%%%%%%%%%%%%%%%%%%%%%%%%%%%%%%%%%%
\section{Refined Estimates for  Degree Growth}
%%%%%%%%%%%%%%%%%%%%%%%%%%%%%%%%%%%%%%%%%%%%%%%%%%%%%%%%%%%%%%%%%%%%%%
For a dominant rational map
\begin{equation}
  \label{eqn:fXdashXHinDivX}
  f : X \longdashrightarrow X
  \quad\text{and ample}\quad H\in\Div(X),
\end{equation}
the definition of $\d(f)$ is equivalent to
\[
  \log\deg_H(f^n) = n\log\d(f) + o(n)\quad\text{as $n\to\infty$.}
\]

\begin{question}
For which $f$ and~$H$ as in~\eqref{eqn:fXdashXHinDivX} is it true that
\[
  \log\deg_H(f^n) = n\log\d(f) + O(n^\e)
  \quad\text{for some $0<\e<1$?}
\]
\end{question}

\begin{question}
For which $f$ and~$H$ as in~\eqref{eqn:fXdashXHinDivX} is it true that
\[
  \log\deg_H(f^n) = n\log\d(f) + O(\log n)?
\]
\end{question}

And for those who want the stars and the moon:

\begin{question}
For which $f$ and~$H$ as in~\eqref{eqn:fXdashXHinDivX} is it true that
the limit
\[
  \lim_{n\to\infty} \frac{\deg_H(f^n)}{\d(f)^n\cdot n^{\ell(f)}}
  \quad\text{exist for some $\ell(f)\in\ZZ_{\ge0}$?}
\]
\end{question}

%%%%%%%%%%%%%%%%%%%%%%%%%%%%%%%%%%%%%%%%%%%%%%%%%%%%%%%%%%%%%%%%%%%%%%
\section{The Arithmetic Degree of an Orbit}
%%%%%%%%%%%%%%%%%%%%%%%%%%%%%%%%%%%%%%%%%%%%%%%%%%%%%%%%%%%%%%%%%%%%%%
For the remainder of this article, we set the convention that
\begin{equation}
  \label{eqn:XQbarfXdashX}
  \left(
  \begin{tabular}{@{}l@{}}
    $X/\Qbar$ is a smooth irreducible projective variety, \\
    $f:X\dashrightarrow X$ is a dominant rational map  defined over~$\Qbar$.\\
  \end{tabular}
  \right)
\end{equation}
We further let
\[
  X(\Qbar)_f:=\bigl\{Q\in X(\Qbar) : \text{$f^n(Q)$ is defined for all $n\ge1$} \bigr\},
\]
and we fix a logarithmic Weil height function
\[
  h_X : X(\Qbar) \longrightarrow \RR_{\ge1}
\]
relative to an ample divisor on~$X$.

\begin{intuition}
$h_X(P)=\text{\# of bits to describe $P$}$.  
\end{intuition}

\begin{definition}
The \emph{arithmetic degree} of the $f$-orbit of the point
$P\in{X}(\Qbar)_f$ is the quantity
\[
  \a(f,P) := \lim_{n\to\infty} h_X\bigl(f^n(P)\bigr)^{1/n}.
\]
Since it is not known in general that this limit exists, we also define
\[
  \overline\a(f,P) := \limsup_{n\to\infty} h_X\bigl(f^n(P)\bigr)^{1/n}.
\]
\end{definition}

\begin{remark}
The triangle inequality easily gives the existence of a constant~$C_f>0$ such
that
\[
h_X\bigl(f(P)\bigr) \le C_f h_X(P) \quad\text{for all $P\in{X}(\Qbar)_f$.}
\]
From this we see that~$\overline\a(f,P)<\infty$. The following precise upper bound
lies much deeper.
\end{remark}

\begin{theorem}
\label{theorem:KSMaled}
\textup{(Kawaguchi--Silverman \cite{MR3456169}, Matsuzawa \cite{arxiv1606.00598})}
\[
  \overline\a(f,P)\le\d(f).
\]
\end{theorem}

\begin{center}
%%    \textbf{Moral.}\quad
    \begin{tabular}{@{}c@{}}
      \textbf{Moral of}\\ \textbf{Theorem \ref{theorem:KSMaled}:}\\
    \end{tabular}
  \framebox{
    $
    \left(\begin{tabular}{@{}c@{}}
      Arithmetic complex-\\
      ity of an orbit\\
    \end{tabular}\right)
    \le
    \left(\begin{tabular}{@{}c@{}}
      Dynamical complex-\\
      ity of the map\\
    \end{tabular}\right)
    $
  }
\end{center}

\begin{conjecture}
\textup{(Kawaguchi--Silverman \cite{MR3456169,MR3483624})}
 The limit defining $\a(f,P)$ converges.
\end{conjecture}

The convergence is known in many situations, including for morphisms and
for many types of maps of surfaces. 

%%%%%%%%%%%%%%%%%%%%%%%%%%%%%%%%%%%%%%%%%%%%%%%%%%%%%%%%%%%%%%%%%%%%%%
\section{Arithmetic Degree Versus Dynamical Degree}
%%%%%%%%%%%%%%%%%%%%%%%%%%%%%%%%%%%%%%%%%%%%%%%%%%%%%%%%%%%%%%%%%%%%%%
\begin{conjecture}
\label{conj:KSdensity}  
\textup{(Kawaguchi--Silverman Density Conjecture\cite{MR3483624})}
If the orbit
\[
  \Orbit_f(P):=\bigl\{f^n(P):n\ge0\bigr\}
\]
is Zariski dense in $X$, then~$\a(f,P)$ exists and satisfies
\[
  \a(f,P)=\d(f). 
\]
\end{conjecture}

\begin{moral}
\(
\left(\begin{tabular}{@{}c@{}}
  Maximal geometric\\
  complexity of an orbit\\
\end{tabular}\right)
\Longrightarrow
\left(\begin{tabular}{@{}c@{}}
  Maximal arithmetic\\
  complexity of the orbit\\
\end{tabular}\right)
\)
\end{moral}

\begin{update}
The density conjecture (Conjecture~\ref{conj:KSdensity}) is known in
some cases, including the following non-exhaustive list:
\begin{enumerate}
\item
Group endomorphisms (homomorphisms composed with translations) of
semi-abelian varieties (extensions of abelian varieties by algebraic tori)
\cite{MR3456169,MR4092861,arxiv1111.5664,MR3614521}.
%% [Silverman]{arxiv1111.5664},
%%      tori, monic monomial maps
%% [Kawaguchi--Silverman]{MR3456169},
%%      abelian varieties, isogenies
%% [Silverman]{MR3614521},
%%      abelian varieties, arbitrary maps
%% [Matsuzawa--Sano]{MR4092861}.
%%      semi-abelian varieties
\item
Endomorphisms of (not necessarily smooth) projective surfaces
\cite{MR3189467,arxiv1701.04369,arxiv1908.01605}.
%% [Kawaguchi--Silverman (2014)]{MR3189467},
%%      smooth projective surfaces, automorphisms
%% [Matsuzawa--Sano--Shibata (2018)]{arxiv1701.04369},
%%      smooth projective surfaces, endomorphisms
%% [Meng--Zhang (2019, preprint)]{arxiv1908.01605}.
%%      projective surfaces, endomorphisms
\item
Extensions to~$\PP^N$ of regular affine automorphisms of $\AA^N$ \cite{MR3189467}.
%% [Kawaguchi--Silverman~(2014)]{MR3189467}.
\item
Endomorphisms of hyperk\"ahler varieties \cite{MR4259157}.
%% [Lesieutre--Satriano (2021)]{MR4259157}.
%% [Matsuzawa--Meng--Shibata--Zhang (2020)]{arxiv2002.10976}.
\item
  Endomorphisms of degree greater than~$1$ of smooth projective
  threefolds of Kodaira dimension~$0$ \cite{MR4259157}.
%% [Lesieutre--Satriano (2021)]{MR4259157}.
\item
Endomorphisms of normal projective varieties having the property that
\text{$\Pic^0\otimes\QQ=0$} and whose nef cone is generated by
finitely many semi-ample integral divisors \cite{arxiv1902.06072}.
%% [Matsuzawa (2020)]{arxiv1902.06072}.
\item
Smooth projective threefolds having
at least one int-amplified\footnote{A morphism $f:X\to{X}$ is \emph{int-amplified} if there
exists an ample Cartier divisor~$H$ such that~$f^*H-H$ is also
 ample.}
endomorphism, and surjective endomorphisms of smooth rationally connected projective varieties
\cite{arxiv2002.10976}.
%% [Matsuzawa--Meng--Shibata--Zhang (2020)]{arxiv2002.10976}.
\end{enumerate}
\end{update}

%%%%%%%%%%%%%%%%%%%%%%%%%%%%%%%%%%%%%%%%%%%%%%%%%%%%%%%%%%%%%%%%%%%%%%
\section{Canonical Heights for Polarized Morphisms}
%%%%%%%%%%%%%%%%%%%%%%%%%%%%%%%%%%%%%%%%%%%%%%%%%%%%%%%%%%%%%%%%%%%%%%

\begin{definition}
Let $f:X\to X$ be a morphism, let $D\in\Div(X)\otimes\RR$ be a
divisor, and suppose that\footnote{We note that~$\d$ need not be an
integer; see for example~\cite{MR3456169,silverman:K3heights} for
applications involving nef divisors~$D$ whose eigenvalues~$\d$ are
algebraic units.}
\[
  f^*D \sim \d D\quad\text{for some $\d>1$.}
\]
We then say that~$D$ is an \emph{eigendivisor} for~$f$,
and we define the associated \emph{canonical height} of a point
\text{$P\in{X}(\Qbar)$} to be the quantity~\cite{callsilv:htonvariety}
\[
\hhat_{f,D}(P) := \lim_{n\to\infty} \frac{1}{\d^n}h_D\bigl(f^n(P)\bigr).
\]
If the divisor~$D$ is ample, then we
say that~$(X,f,D)$ is a \emph{polarized dynamical system}.
\end{definition}

We also recall the standard notation 
\begin{align}
  \label{eqn:PrePerf}
  \PrePer(f) &:= \{\text{preperiodic points of $f$}\} \notag\\
  &= \bigl\{ P\in X : \text{$\Orbit_f(P)$ is finite} \bigr\}.
\end{align}

Standard properties of $\hhat_{f,D}$ include:\footnote{For
the third property, we are using the fact~\eqref{eqn:XQbarfXdashX}
that we are working over~$\Qbar$, although it holds more generally
for fields and heights having the Northcott property.}
\begin{align*}
  \hhat_{f,D}(P)&=h_D(P)+O(1); \\
  \hhat_{f,D}\bigl(f(P)\bigr)&=\d\hhat_{f,D}(P); \\
  \hhat_{f,D}(P)&=0\Longleftrightarrow P\in\PrePer(f) \quad\text{for $D$ ample.}
\end{align*}

\begin{conjecture}
  \textup{(Dynamical Lehmer Conjecture)}
Let  $D$ be ample. Then there exists a constant $C(f,D)>0$ so that
\[
\hhat_{f,D}(P) \ge \frac{C(f,D)}{\bigl[ \QQ(P) : \QQ \bigr]}
\quad\text{for all $P\in{X}(\Qbar)\setminus\PrePer(f)$.}
\]
\end{conjecture}

%%%%%%%%%%%%%%%%%%%%%%%%%%%%%%%%%%%%%%%%%%%%%%%%%%%%%%%%%%%%%%%%%%%%%%
\section{Shibata's Ample Canonical Height}
%%%%%%%%%%%%%%%%%%%%%%%%%%%%%%%%%%%%%%%%%%%%%%%%%%%%%%%%%%%%%%%%%%%%%%
Let $f:X\to{X}$ be a dominant morphism with $\d(f)>1$, let
\[
h_X:X(\Qbar)\to\RR_{\ge1}
\]
be a logarithmic Weil height relative to an ample divisor.
We define~$\ell(f)$ to be the quantity
\[
\ell(f) = \inf\left\{
\ell\ge0 : \sup_{n\ge1}\;\; \frac{h_X\bigl(f^n(P)\bigr)}{n^{\ell}\cdot\d(f)^n} < \infty \right\},
\]
where we formally define the infimum of the empty set to be~$\infty$.

\begin{definition}
The (\emph{lower}) \emph{ample canonical height} is 
\[
\hhatinf_f:X(\Qbar)\to\RR_{\ge0},\quad
\hhatinf_f(P) := \liminf_{n\to\infty} \frac{h_X\bigl(f^n(P)\bigr)}{n^{\ell(f)}\cdot\d(f)^n},
\]
where we formally set~$\hhatinf_f(P)=0$ if~$\ell(f)=\infty$.  
\end{definition}

\begin{conjecture}
\label{conjecture:shibata}
\textup{(Shibata \cite{arxiv1710.05278})}
Let $f:X\to{X}$ be a dominant morphism with $\d(f)>1$.
\begin{parts}
\Part{(a)}
The quantity~$\ell(f)$ is a non-negative integer.
\Part{(b)}
For every number field~$K/\QQ$ over which~$X$ and~$f$ are defined, the set
\begin{equation}
  \label{eqn:PXKhinf}
 \bigl\{ P\in X(K) : \hhatinf_f(P) = 0 \bigr\} 
\end{equation}
is not Zariski dense in $X$.
\end{parts}
%% \textup(N.B. For this conjecture, the map~$f$ is a morphism.\textup)
\end{conjecture}

\begin{remark}
The set \eqref{eqn:PXKhinf} is independent of the choice of ample height~$h_X$.
\end{remark}

%%%%%%%%%%%%%%%%%%%%%%%%%%%%%%%%%%%%%%%%%%%%%%%%%%%%%%%%%%%%%%%%%%%%%%
\section{Shibata's Conjecture Implies the Kawaguchi--Silverman Density Conjecture}
%%%%%%%%%%%%%%%%%%%%%%%%%%%%%%%%%%%%%%%%%%%%%%%%%%%%%%%%%%%%%%%%%%%%%%
Suppose that~$\a(f,P)$ exists and
satisfies~\text{$\a(f,P)<\d(f)$}. For notational convenience, we let
\[
\e(f,P) := \d(f) - \a(f,P) > 0.
\]
Then for sufficiently large~$n$,
\[
\begin{aligned}
  h_X\bigl(f^n(P)\bigr)^{1/n}
  &\le \a(f,P) + \tfrac12 \e(f,P) \\
  &= \d(f) - \tfrac12 \e(f,P).
\end{aligned}
\]
Hence
\[
\begin{aligned}
  \hhatinf_f(P)
  &= \liminf_{n\to\infty} \frac{h_X\bigl(f^n(P)\bigr)}{n^{\ell(f)}\cdot\d(f)^n} \\
  &\le \liminf_{n\to\infty}
  \frac{ \bigl( \d(f) - \frac12\e(f,P) \bigr)^n}{n^{\ell(f)}\cdot\d(f)^n} \\
  &= \liminf_{n\to\infty} \frac{1}{n^{\ell(f)}} \left(1 - \frac{\e(f,P)}{2\d(f)} \right)^n \\
  &= 0.
\end{aligned}
\]
Since $\a\bigl(f,f^n(P)\bigr)=\a(f,P)$, we see that a strict
inequality $\a(f,P)<\d(f)$ implies that
\[
 \Orbit_f(P) \subseteq \bigl\{ Q\in X(K) : \hhatinf_f(Q)=0 \bigr\}.
\]
Therefore Shibata's conjecture (Conjecture~\ref{conjecture:shibata}),
which says that this last set is not Zariski dense, implies the
Kawaguchi--Silverman density conjecture
(Conjecture~\ref{conj:KSdensity}) for morphisms.

%%%%%%%%%%%%%%%%%%%%%%%%%%%%%%%%%%%%%%%%%%%%%%%%%%%%%%%%%%%%%%%%%%%%%%
\section{Other Types of Growth Rates for $h_X\bigl(f^n(P)\bigr)$?}
%%%%%%%%%%%%%%%%%%%%%%%%%%%%%%%%%%%%%%%%%%%%%%%%%%%%%%%%%%%%%%%%%%%%%%
The map
\[
  f:\PP^2\dashrightarrow\PP^2,\quad f(x,y,z)=[xy+xz,yz+z^2,z^2]
\]
is interesting. It satisfies
\[
  \deg(f^n)=n+1,
\]
so $\d(f)=1$, and the orbit of the point $P=[1,0,1]$ satisfies
\[
  f^n(P) = [n!,n,1],
\]
so
\[
  h\bigl(f^n(P)\bigr)=\log(n!)\sim n\log n.
\]

\begin{question}
For dominant rational maps $f:\PP^N\dashrightarrow\PP^N$,
what types of growth are possible? For example
is it possible to have:
\par
\begin{tabular}{ll}
  (a)& $h\bigl(f^n(P)\bigr)\sim n^i(\log n)^j$ for some $j\ge2$? \\[2\jot]
  (b)& $h\bigl(f^n(P)\bigr)\sim \d(f)^n n^i(\log n)^j$ for some $j\ge1$ and $\d(f)>1$? \\
\end{tabular}
\end{question}

%%%%%%%%%%%%%%%%%%%%%%%%%%%%%%%%%%%%%%%%%%%%%%%%%%%%%%%%%%%%%%%%%%%%%%
\section{Canonical Heights and Heights of Dynamical Systems}
\label{section:moduliht}
%%%%%%%%%%%%%%%%%%%%%%%%%%%%%%%%%%%%%%%%%%%%%%%%%%%%%%%%%%%%%%%%%%%%%%
Let $\Moduli_d^N:=\End_d(\PP^N)/\!/\PGL_{N+1}$ denote
the moduli space of degree~$d$ dynamical systems (morphisms)~$\PP^N\to\PP^N$,
and fix an ample height
\begin{equation}
  \label{eqn:hmoduli}
  h_\Moduli : \Moduli_d^N(\Qbar) \longrightarrow \RR_{\ge1}.
\end{equation}

\begin{conjecture}
\textup{(Dynamical Lang Height Conjecture)}
Let $K/\QQ$ be a number field. There are constants $C_1(K,N,d)>0$ and
$C_2(K,N,d)$ so that for all $f\in\Moduli_d^N(K)$ and all
$P\in\PP^N(K)$ with Zariski dense $f$-orbit,
\[
  \hhat_f(P) \ge C_1h_\Moduli(f) - C_2.
\]
\end{conjecture}

%%%%%%%%%%%%%%%%%%%%%%%%%%%%%%%%%%%%%%%%%%%%%%%%%%%%%%%%%%%%%%%%%%%%%%
\section{Critical Heights and Heights of Dynamical Systems}
%%%%%%%%%%%%%%%%%%%%%%%%%%%%%%%%%%%%%%%%%%%%%%%%%%%%%%%%%%%%%%%%%%%%%%

For this section we restrict to~$\PP^1$.
We recall from~\eqref{eqn:PrePerf}
the notation~$\PrePer(f)$ for the preperiodic points of~$f$,
and we let~$\Crit_f$ denote the critical points of~$f$.

\begin{definition}
The \emph{critical height} of $f\in\Moduli_d^1(\Qbar)$ is
\[
  \hhat^\crit(f) := \sum_{P\in\Crit_f} \hhat_f(P).
\]
\end{definition}

Notice that we have
\[
  \hhat^\crit(f)=0
  \Longleftrightarrow
  \Crit_f\subset\PrePer(f)
  \Longleftrightarrow
  \text{$f$ is PCF}.
\]

\begin{definition}
\label{definition:commensurate}
We say that non-negative
real-valued functions~$\f$ and~$\psi$ are \emph{commensurate}, and we write~$\f\asymp\psi$, if there
are positive constants~$c_1,c_2,c_3,c_4$ such
that
\[
c_1\f(x)-c_2\le\psi(x)\le{c_3}\f(x)+c_4\quad\text{for all~$x$.}
\]
\end{definition}

\begin{theorem}
\label{thm:ing}
\textup{(Ingram \cite{MR3799700})}
The moduli height and the critical height 
are commensurate \textup(Definition~\ref{definition:commensurate}\textup)
away from the locus of Latt{\`e}s maps,
\[
\hhat^\crit  \asymp h_\Moduli
\quad\text{on the set of non-Latt{\`e}s maps in~$\Moduli_d^1(\Qbar)$.}
\]
\end{theorem}

\begin{question}
How might we generalize Theorem~\ref{thm:ing} to
$\Moduli_d^N$? What replaces $\hhat^\crit$?
We shall give one possible answer in Section~\ref{section:crithtPN}.
\end{question}

%%%%%%%%%%%%%%%%%%%%%%%%%%%%%%%%%%%%%%%%%%%%%%%%%%%%%%%%%%%%%%%%%%%%%%
\section{Higher Order Dynamical Degrees} %% ???
\label{section:higherorderdyndeg}
%%%%%%%%%%%%%%%%%%%%%%%%%%%%%%%%%%%%%%%%%%%%%%%%%%%%%%%%%%%%%%%%%%%%%%
Let~$X$ be a smooth irreducible projective variety of dimension~$N$,
let \text{$f:X\dashrightarrow{X}$} be a dominant rational map, and let~$H$ be
an ample divisor.

\begin{definition}
The \emph{$k$th dynamical degree of ${f}$} is
\begin{equation}
\label{eqn:dkf}
\d_k(f) := \lim_{n\to\infty} \Bigl( (f^n)^*(H^k) \cdot H^{N-k} \Bigr)^{1/n}.
\end{equation}
\end{definition}

See~\cite{arXiv:1701.07760,MR2119243,MR4048444,MR3449182}
for proofs that the limit~\eqref{eqn:dkf} exists and has various
desirable properties.

\begin{intuition}
The $k$th dynamical degree~$\d_k(f)$ measures the dynamical complexity
of~$f$ in codimension~$k$.  For example,
\begin{equation}
  \label{eqn:dnfeqtopdeg}
  \d_N(f) = \text{topological degree of $f$} = \#f^{-1}(\text{generic point}).
\end{equation}
\end{intuition}

\begin{theorem}
\label{thm:gue}
\textup{(Guedj \cite{MR2179389})}
Dynamical degrees form a log concave sequence, i.e., $\d_{i-1}\d_{i+1}\le\d_i^2$.
In particular, for some $k$ we have
\[
\d_1(f) \le \d_2(f) \le \cdots \le \d_k(f)
\quad\text{and}\quad
\d_k(f) \ge \d_{k+1}(f) \ge\cdots\ge \d_N(f).
\]
\end{theorem}

For a description of the full sequence of dynamical degrees when~$f$
is a monomial map, see \cite[Favre--Wulcan~(2012)]{arxiv1011.2854} and
\cite[Lin~(2012)]{arxiv1010.6285}.

\begin{question}
\label{question:dkfalgint}
Are all $\d_k(f)$ algebraic integers?  
\end{question}

\begin{update}
We have~$\d_N(f)\in\ZZ$ from~\eqref{eqn:dnfeqtopdeg}, but the
existence of maps with transcendental~$\d_1(f)$ as proven
in~\cite{belldillerjonssonkrieger2021,arxiv1907.00675} shows that
Question~\ref{question:dkfalgint} has a negative answer
for~$k=1$. This leave Question~\ref{question:dkfalgint} open
for~\text{$2\le{k}<N$}.
\end{update}

\begin{question}
Is there any way to generalize Theorem~\ref{thm:gue} to some sort of
higher order arithmetic degrees? We consider this question in the next section.
\end{question}

%%%%%%%%%%%%%%%%%%%%%%%%%%%%%%%%%%%%%%%%%%%%%%%%%%%%%%%%%%%%%%%%%%%%%%
\section{Higher Order Arithmetic Degree}
%%%%%%%%%%%%%%%%%%%%%%%%%%%%%%%%%%%%%%%%%%%%%%%%%%%%%%%%%%%%%%%%%%%%%%
Let~$X/\Qbar$ be a smooth irreducible projective variety of dimension~$N$, and let
$f:X\to{X}$ be a morphism. For~$k\ge2$, the ``natural'' way to define
the~$k$th arithmetic degree~$\a_k(f,P)$ of a point~$P$ using
Arakelov-type intersection theory (conjecturally) yields
\[
  \a_k(f,P)=1\quad\text{for all $P$,}
\]
which is not very interesting.  Indeed, in the definition
of~$\a_1(f,P)$, the point~$P$ has dimension~$1$ and a divisor~$H$ has
codimension~$1$, so we expect the arithmetic intersection of~$f^n(P)$
and~$H$ to be quite large in general.  But if we replace~$H$ with,
say,~$H_1\cdot{H_2}$, then the intersection of $H_1\cdot{H_2}$
with~$f^n(P)$ is likely to be small for all $n\ge0$.

One possible solution is to replace the point~$P$ with a higher dimensional subvariety.
So we want to define the $k$th arithmetic degree to be a function
\[
%% \a_k(f,\,\cdot\,) : \{ \text{$k-1$ dim'l irreducible subvarieties} \} \longrightarrow \RR_{\ge0}
\a_k(f,\,\cdot\,) : \left\{
\begin{tabular}{@{}l@{}}
$k-1$ dimensional irredu-\\ cible subvarieties of $X$\\
\end{tabular}
\right\} \longrightarrow \RR_{\ge0}
\]
whose domain is the set of irreducible subvarieties of~$X$ of geometric
dimension~\text{$k-1$}.

There is a theory that assigns a height to each subvariety
\[
  Z \subseteq X,
\]
in particular when~$X=\PP^N$. Indeed, there are several formulations,
including a height for $Z\subset\PP^N$ using the Chow coordinates
of~$Z$, and a height $h_{X,\overline{L}}(Z)$ relative to a metrized
line bundles~$\overline{L}$ via Arakelov theory due to Faltings,
Zhang, and Bost--Gillet--Soul{\'e}; see for example~\cite{MR1260106}.

\begin{definition}
Let $f:X\to{X}$ be a morphism,  let $Z\subset{X}$ be a \text{$k-1$}-dimensional
irreducible subvariety of~$X$, and let~$\overline{L}$ be a metrized line bundle on~$X$
as defined for example in~\cite{MR1260106}.
The \emph{arithmetic degree} of $Z\subseteq{X}$ is
\[
  \a_k(f,Z) := \lim_{n\to\infty} h_{X,\overline{L}}\bigl(f^n(Z)\bigr)^{1/n}.
  \]
\end{definition}  

\begin{question}
Does the limit $\a_k(f,Z)$ converge?  
\end{question}

\begin{question}
Is there a natural upper bound for $\a_k(f,Z)$ in terms of
$\d_1(f),\ldots,\d_{1+\dim Z}(f)$? If so, when is this upper bound
attained?
\end{question}

\begin{example}
(Kawaguchi--Silverman, unpublished)
Let
\[
f:\PP^N\longdashrightarrow\PP^N
\]
be a dominant {monomial} map, i.e., a map whose coordinate functions
are monomials, and let $Z\subset\PP^N$ be an irreducible hypersurface
that is not a coordinate hyperplane. Then
\[
\overline\a_N(f,Z) \le \min \bigl\{ \d_{N-1}(f), \d_N(f) \bigr\}.
\]
Further, there are examples with $N=2$ satisfying:
\begin{align*}
(1)~\a_N(f,Z) &= \d_1(f) < \d_2(f); \\
(2)~\a_N(f,Z) &= \d_2(f) < \d_1(f).
\end{align*}
The proofs are elementary triangle inequality computations.
\end{example}

%%%%%%%%%%%%%%%%%%%%%%%%%%%%%%%%%%%%%%%%%%%%%%%%%%%%%%%%%%%%%%%%%%%%%%
\section{Height Lower Bounds and the Bogomolov Property}
%%%%%%%%%%%%%%%%%%%%%%%%%%%%%%%%%%%%%%%%%%%%%%%%%%%%%%%%%%%%%%%%%%%%%%
We fix a polarized dynamical system $(X,f,D)$, i.e.,
\[
  f:X\to X,\quad \text{$D$ ample},\quad
  \text{$f^*D\sim\d{D}$ for some $\d>1$.}
\]
For a subvariety $Z\subseteq{X}$ and any $\e>0$, we let
\[
  \overline{Z_{f,D}(\e)} := \text{Zariski closure of }
  \bigl\{ P\in Z(\Qbar): \hhat_{X,f,D}(P) < \e \bigr\}. 
\]
  
\begin{definition}
\label{definition:bogomolovproperty}
A subvariety $Z\subseteq{X}$ has the \emph{Bogomolov property}
(relative to $f$ and $D$) if there is an $\e>0$ such that
\[
  \overline{Z_{f,D}(\e)} \ne Z.
\]
\end{definition}

\begin{example}
Let $X$ be an abelian variety, and assume that $Z\subset{X}$ is not a
translate of an abelian subvariety by a torsion point. Then
Zhang~\cite{MR1609518}, Ullmo~\cite{MR1609514}, and
David--Philippon~\cite{MR1478502} prove that~$Z$ has the Bogomolov
property.
\end{example}

\begin{example}
Let $X=(\PP^1)^N$, let $f$ be a dominant endomorphism, and let
$Z\subset{X}$ be an irreducible subvariety that is
not~$f$-preperiodic. Then Ghioca--Nguyen--Ye prove~\cite{MR3826461}
that~$Z$ has the Bogomolov property.
\end{example}

\begin{remark}
As shown by a construction of Ghioca, Tucker, and Zhang~\cite[Remark~2.3]{MR2854724},
\[
  \overline{Z\cap\PrePer(f)}=Z \quad
  \text{does \textbf{not} imply that $Z$ is $f$-preperiodic.}
\]
In other words, the naive statement of the dynamical Manin--Mumford conjecture is not true.
See~\cite{MR4205770} for a formulation that may be true.
\end{remark}

%%%%%%%%%%%%%%%%%%%%%%%%%%%%%%%%%%%%%%%%%%%%%%%%%%%%%%%%%%%%%%%%%%%%%%
\section{The Bogomolov Canonical Height of a Subvariety}
%%%%%%%%%%%%%%%%%%%%%%%%%%%%%%%%%%%%%%%%%%%%%%%%%%%%%%%%%%%%%%%%%%%%%%
Definition~\ref{definition:bogomolovproperty} suggests defining the
Bogomolov height of~$Z$ to be the largest~$\e$ for which~$Z$ has the
Bogomolov property.

\begin{definition}
The \emph{Bogomolov height of ${Z}$} (relative to $f$ and $D$) is
\[
 \hhat^\Bcal_{X,f,D}(Z) := \sup_{\emptyset\ne U\subseteq Z} \;\; \inf_{\strut P\in U(\Qbar)} \hhat_{X,f,D}(P),
\]
where $U$ ranges over Zariski open subsets of~$Z$. And if
$Z=\sum{n_iZ_i}$ is a formal sum of equidimensional subvarieties,
we extend linearly,
\[
  \hhat^\Bcal_{X,f,D}(Z)=\sum{n_i}\hhat^\Bcal(Z_i).
\]
\end{definition}
Equivalently, since
\[
\hhat^\Bcal_{X,f,D}(Z) = \sup \Bigl(
\{0\} \cup \bigl\{ \e>0 : \overline{Z_{f,D}(\e)}\ne Z \bigr\} \Bigr),
\]
we see that that
\[
    \text{$Z$ has the Bogomolov property}
    \;\Longleftrightarrow\;
    \hhat^\Bcal_{X,f,D}(Z)>0.
\]

\begin{remark}
Maybe we should say that~$Z$ is \emph{formally preperiodic}
if $\hhat^\Bcal_{X,f,D}(Z)=0$?
\end{remark}

\begin{definition}
Define the \emph{canonical height of~$Z$} to be the quantity  
\[
  \hhat_{X,f,D}(Z) := \lim_{n\to\infty} \frac{1}{\d^n} h_{X,D}\bigl(f^n(Z)\bigr).
\]
\end{definition}

\begin{theorem}
\textup{(Zhang \cite{MR1609518,MR1311351})}
\begin{parts}
\Part{(a)}
The canonical height $\hhat_{X,f,D}(Z)$ converges. \textup(N.B. We
are working with a polarized dynamical system.\textup)
\Part{(b)}
The canonical height and the Bogomolov height are
are commensurate \textup(Definition~\ref{definition:commensurate}\textup),
\[
\hhat_{X,f,D} \asymp \hhat_{X,f,D}^\Bcal
\quad\text{on the set of irreducible subvarieties $Z\subset X$.}
\]
\end{parts}
\end{theorem}

%%%%%%%%%%%%%%%%%%%%%%%%%%%%%%%%%%%%%%%%%%%%%%%%%%%%%%%%%%%%%%%%%%%%%%
\section{The Critical Height of an Endomorphism of $\PP^N$}
\label{section:crithtPN}
%%%%%%%%%%%%%%%%%%%%%%%%%%%%%%%%%%%%%%%%%%%%%%%%%%%%%%%%%%%%%%%%%%%%%%
Let $f:\PP^N\to\PP^N$ be an endomorphism, and let
\[
  \Crit_f := \text{(the critical locus of $f$)} \in \Div(\PP^N),
\]
where we view~$\Crit_f$ as an effective divisor with appropriate multiplicities.

\begin{definition}
The \emph{critical height of ${f:\PP^N\to\PP^N}$} is
the quantity  
\[
  \hhat^\crit(f) := \hhat^\Bcal_{X,f,\Ocal(1)}(\Crit_f).
\]
\end{definition}

This gives a well-defined function
\[
  \hhat^\crit : \Moduli_d^N(\Qbar) \longrightarrow \RR_{\ge0},
\]
and we might say that~$f$ is ``formally PCF'' if $\hhat^\crit(f)=0$.

\begin{question}
Is it possible to have $\hhat^\crit(f)=0$ with $\Crit_f$ not being
an $f$-preperiodic subvariety?  
\end{question}

We recall from Section~\ref{section:moduliht} that~$h_\Moduli$ denotes
an ample height on the moduli space~$\Moduli_d^N(\Qbar)$ of degree~$d$
endomorphisms of~$\PP^N$.  The following conjecture generalizes
Ingram's theorem (Theorem~\ref{thm:ing}) to higher dimension.

\begin{conjecture}
There is a Zariski closed set $L_d^N\subsetneq\Moduli_d^N$ so that the
critical height and the moduli height are commensuarate
\textup(Definition~\ref{definition:commensurate}\textup) on the
complement of~$L_d^N$,
\[
\hhat^\crit  \asymp h_\Moduli
\quad\text{on the set of maps $f\in\Moduli_d^N(\Qbar)\setminus L_d^N$.}
\]
\end{conjecture} 

\begin{remark}
The upper bound $\hhat^\crit\le{c_1}h_\Moduli+c_2$ is probably not too hard;
cf.\ the proof of~\cite[Theorem~6.31]{MR2884382}. The lower bound, if
true, seems considerably more difficult, in particular because we have
no clear understanding of the structure of the exceptional set~$L_d^N$
of higher dimensional Latt{\`e}s-type maps.
\end{remark}

%%%%%%%%%%%%%%%%%%%%%%%%%%%%%%%%%%%%%%%%%%%%%%%%%%%%%%%%%%%%%%%%%%%%%%
\section{Height relative to subvarieties of codimension at least two}
%%%%%%%%%%%%%%%%%%%%%%%%%%%%%%%%%%%%%%%%%%%%%%%%%%%%%%%%%%%%%%%%%%%%%%

\begin{definition}
Let $X/\Qbar$ be a smooth projective variety, let~$Z\subset{X}$ be a
subvariety (or more generally, a closed subscheme). Then a Weil height
on $X$ relative to~$Z$ is a function
\[
h_{X,Z}: X(\Qbar)\setminus Z \longrightarrow \RR
\]
defined by fixing a blow-up $\pi:\hat{X}\to{X}$ of~$X$ along~$Z$,
letting~$\hat{Z}$ be the exceptional divisor on~$\hat{X}$, choosing a
Weil height~$h_{\hat{X},\hat{Z}}$ on~$\hat{X}(\Qbar)$ associated
to~$\hat{Z}\in\Div(\hat{X})$, and setting
\[
h_{X,Z}(P) = h_{\hat{X},\hat{Z}}(\hat{P}) \quad\text{where $\hat{P}:=\pi^{-1}(P)$.}
\]
See~\cite[Theorems~2.1 and~4.1]{silverman:arithdistfunctions} for
details of this construction.\footnote{\emph{Errata
for~\cite{silverman:arithdistfunctions}\/}: Item~(iii) on page~196
should be~$\mathcal{I}_{X\cup{Y}}$ instead
of~$\mathcal{I}_{X\cap{Y}}$.  In item~(e) on pages~197 and~203, the
subscript on the local and global heights in the middle of the
inqualities should be~$X\cup{Y}$ instead
of~$X\cap{Y}$. Proposition~6.2 is incorrect; see~\cite{MR4247874} for
a corrected version.} The intuition behind this definition is
that~$h_{X,Z}(P)$ is a sum of local heights $\lambda_{X,Z,v}(P)$ that
are supposed to measure the~$v$-adic distances from~$P$ to~$Z$. Since
we already know how to define local heights associated to divisors, an
easy way to do this is to replace~$Z$ with the divisor~$\hat{Z}$ and
use the classical Weil construction.
\end{definition}

The following question was inspired by a question of Matsuzawa in
which~$f$ was a rational map and one studied the height~$h_{X,I(f)}$
relative to the indeterminacy locus~$I(f)$ of~$f$.

\begin{definition}
Let $f:X\dashrightarrow X$ be a dominant rational map, let
$Y\subseteq{X}$ be a subvariety, and let~$P\in{X}$. We say that~$Y$ is
an \emph{accumulating subvariety for the $f$-orbit of~$P$} if
\[
  \#\bigl\{ n\in\ZZ_{\ge0} : f^n(P)\in Y \bigr\} = \infty.
\]
\end{definition}

\begin{question}
\label{ques:hXZvshXH}
Let $X/\Qbar$ be a smooth projective variety, let~$Z\subset{X}$ be an
irreducible subvariety of codimension at least~$2$, let $f:X\to{X}$ be
an endomorphism defined over~$\Qbar$, and let $P\in{X}(\Qbar)$ be a
point whose~$f$-orbit has no non-trivial accumulating subvarieties.
We ask for conditions on~$X$,~$Z$,~$f$ which guarantee that
\begin{equation}
  \label{eqn:hZhh0}
  \lim_{n\to\infty} \frac{h_{X,Z}\bigl(f^n(P)\bigr)}{h_{X,H}\bigl(f^n(P)\bigr)}=0
  \quad\text{for all $P\in X(\Qbar)$.}
\end{equation}
We pose two specific questions.
\begin{parts}
\Part{(a)}
Is~\eqref{eqn:hZhh0} true if $Z\cap\Crit_f=\emptyset$?
\Part{(b)}
Is~\eqref{eqn:hZhh0} true if $Z\cap\Crit_f$ contains no~$f$-periodic points?
\end{parts}
\end{question}

\begin{example}
Here is an example where the limit~\eqref{eqn:hZhh0} in
Question~\ref{ques:hXZvshXH} is true. Let~$a,b\in\ZZ$ be non-zero
integers that are multiplicatively independent. Consider the following
map, subvariety, and point on~$\PP^2$:
\[
  f(x,y,z)=[ax,by,z],\quad Z=\bigl\{[1,1,1]\bigr\},\quad P=[1,1,1].
\]
It is not hard to check that if $a,b\in\ZZ$ with $(a,b)\ne(1,1)$, then~\cite[Example in Section~4]{silverman:arithdistfunctions}
\[
h_Z(a,b,1) = \log\gcd(a-1,b-1).
\]
Since
\[
f^n(P) = [a^n,b^n,1],
\]
we see that
\begin{align*}
h_{\PP^2,Z}\bigl(f^n(P)\bigr) &= \log\gcd(a^n-1,b^n-1),\\
h_{\PP^2,H}\bigl(f^n(P)\bigr) &= \log\max(|a^n-1|,|b^n-1|) \sim n \log\max\{|a|,|b|\}.
\end{align*}
Hence~\eqref{eqn:hZhh0} becomes the statement that
\begin{equation}
  \label{eqn:limanbn}
\lim_{n\to\infty} \frac{\log\gcd(a^n-1,b^n-1)}{n} = 0,
\end{equation}
which is a theorem of Bugeaud, Corvaja, and
Zannier~\cite{MR1953049}. However, it should be noted that the proof of~\eqref{eqn:limanbn}
uses Schmidt's subspace theorem, so proving the elementary
looking inequality~\eqref{eqn:limanbn} currently requires some highly
non-elementary machinery.  See~\cite[Section~2]{MR2162351} for further
details.
\end{example}

\begin{example}
And here is an observation of Chong Gyu (Joey) Lee (private
communication) showing that some condition on~$f$,~$P$, or~$Z$ is
required for the validity of the limit~\eqref{eqn:hZhh0}. Fix $d\ge2$
and $a,b\in\ZZ_{\ne0}$, and let
\[
f(x,y,z)=[x^d,y^d,z^d],\quad Z=\bigl\{[0,0,1]\bigr\},\quad P=[a,b,1].
\]
Then
\[
h_Z\bigl(f^n(P)\bigr) = h_Z\bigl(a^{d^n},b^{d^n},1\bigr)
= \log\gcd(a^{d^n},b^{d^n}),
\]
so
\[
\frac{h_Z\bigl(f^n(P)\bigr)}{h_H\bigl(f^n(P)\bigr)}
= \frac{\log\gcd(a^{d^n},b^{d^n})}{\log\max(|a^{d^n}|,|b^{d^n}|)}
= \frac{\log\gcd(a,b)}{\log\max(|a|,|b|)}.
\]
Hence the limit~\eqref{eqn:hZhh0} is non-zero if $\gcd(a,b)\ge2$. Of
course, in this example, the set~$Z$ consists of a critical fixed
point of~$f$.
\end{example}

\begin{proposition}
In the setting of Question~$\ref{ques:hXZvshXH}$, assume that~$X$ has
trivial canonical divisor.  Then Vojta's conjecture~\cite{MR883451}
implies that the limit formula~\eqref{eqn:hZhh0} is true.
\end{proposition}
\begin{proof}
We let~$\hat{X}\to{X}$ be the blow-up of~$X$ along~$Z$, and we
let~$\hat{Z}\in\Div(\hat{X})$ be the proper transform of~$Z$.  The
canonical divisor on the blow-up is given in general by
\[
K_{\hat{X}} =\pi^* K_X + \hat{Z},
\]
so our assumption that~$K_X=0$ gives
\[
K_{\hat{X}} = \hat{Z}.
\]
\par
Vojta's conjecture on~$\hat{X}$ with
the divisor~$D=0$ says that for every $\e>0$ there
is a Zariski closed subset $E(\e)\subsetneq{\hat{X}}$ such that
\[
h_{\hat{X},K_{\hat{X}}}(\hat{Q}) \le \e h_{\hat{X},\hat{H}}(\hat{Q})
\quad\text{for all $\hat{Q}\in{\hat{X}(K)}\setminus{E}(\e)$.}
\]
Hence for $Q\in{X}(K)\setminus(Z\cup\pi(E(\e)))$, we 
let~$\hat{Q}:=\pi^{-1}(Q)$ and compute
\[
h_{X,Z}(Q) = h_{\hat X,\hat Z}(\hat Q) = h_{\hat X,K_{\hat X}}(\hat Q)
\le \e h_{\hat{X},\hat{H}}(\hat{Q}) = \e h_{X,H}(Q).
\]
It follows that
\[
  \limsup_{\substack{Q \in X(K) \setminus (Z\cup E(\e))}}
    \frac{h_{X,Z}(Q)}{h_{X,H}(Q)} \le \e.
\]
We now set~$Q=f^n(P)$ and use the assumption that~$\Orbit_f(P)$ has
no accumulating subvarieties, which means in particular that there
are only finitely many~$n$ such that $f^n(P)\in{Z}\cup{E(\e)}$,
to conclude that
\[
  \limsup_{n\to\infty}
    \frac{h_{X,Z}\bigl(f^n(P)\bigr)}{h_{X,H}\bigl(f^n(P)\bigr)} \le \e.
\]
This is true for every~$\e>0$, which completes the proof
that the limit exists~\eqref{eqn:hZhh0} and is equal to~$0$.
\end{proof}

%%%%%%%%%%%%%%%%%%%%%%%%%%%%%%%%%%%%%%%%%%%%%%%%%%%%%%%%%%%%%%%%%%%%%%%%
% Acknowldegements and Bibliography
%%%%%%%%%%%%%%%%%%%%%%%%%%%%%%%%%%%%%%%%%%%%%%%%%%%%%%%%%%%%%%%%%%%%%%%%

\begin{acknowledgement}
I would like to thank the Simons Foundation for providing the lovely
setting at which this talk was presented, the symposium organizers
Laura DeMarco and Mattias Jonsson for inviting me to speak, and the
referee for their helpful comments and suggestions.
\end{acknowledgement}

%% \begin{thebibliography}{99}
%% \itemsep=\smallskipamount
%% \end{thebibliography}

%% \bibliographystyle{plain}
%% \bibliography{/Users/joesilverman/Dropbox/AAJHS/Book/ADS/ArithDyn}

\end{document}